\documentclass[12,reqno]{amsart}%
\usepackage{amsfonts}
\usepackage{amsthm}
\usepackage[centertags]{amsmath}
\usepackage{amssymb}
\usepackage{enumerate}
\usepackage{graphicx}
\usepackage[bookmarks=true, pdfstartview={FitH}, pdftitle={Paper}, pdfauthor={Sajid Iqbal}]{hyperref}

\providecommand{\U}[1]{\protect\rule{.1in}{.1in}}
\theoremstyle{plain}
\newtheorem{theorem}{Theorem}[section]
\newtheorem{corollary}[theorem]{Corollary}

\newtheorem{example}[theorem]{Example}

\newtheorem{remark}[theorem]{Remark}

\numberwithin{equation}{section}

\def\proof{\noindent {\it Proof.} }
\def\qed{\hfill $\blacksquare$}

\def\proof{\noindent {\it Proof.} }
\def\qed{\hfill $\blacksquare$}

\begin{document}
\title[weighted opial inequalities for widder derivatives]{weighted opial inequalities for widder derivatives}
\author[]{}
\address{}
\email{}
\author[S. Iqbal$^{1}$]{Sajid Iqbal$^{1}$}
\address{1-Department of Mathematics, University of Sargodha(Sub-Campus Bhakkar), Bhakkar, Pakistan}
\email{sajid$\_$uos2000@yahoo.com}

\author[J. Pe\v cari\'c$^{2}$]{Josip Pe\v cari\'c$^{2}$}
\address{2-Faculty of Textile Technology, University of Zagreb, Prilaz baruna Filipovi\'{c}a 28a, 10000 Zagreb, Croatia}
\email{pecaric@element.hr}
\date{}
\dedicatory{}
\begin{abstract}
In this paper, we establish some new general Opial inequalities for Widder derivatives.
\end{abstract}
\maketitle



\section{Introduction}
The Opial inequality, which appeared in \cite{OP}, is of great interest in differential and difference equations and other areas of mathematics, and has attracted a great deal of attention in the recent literature (see, for instance, \cite{AG}, \cite{AP1}, \cite{AP2}, \cite{AN}, \cite{ANAPE}, \cite{CDP}, \cite{PACH}). Recall that the original inequality \cite{OP} (see also \cite[p. 114]{MPF}) states the following:
\begin{theorem}
Let $a > 0$. If $f \in C^{1}[0,a]$ with $f(0) = f(a) = 0 $ and $f(t) > 0$ on $(0,a)$, then
\begin{equation}
\int\limits_{0}^{a}|f(t)f^{'}(t)|dt\leq\frac{a}{4}\int\limits_{0}^{a}(f^{'}(t))^{2}dt.
\end{equation}
The constant $\frac{a}{4}$ is the best possible.
\end{theorem}

In \cite{ANPA}, Anastassiou established some various $L_{p}$ form of Opial type inequalities for Widder derivatives. Anastassiou and Pe\v cari\'c in \cite{ANAPE} gave the very general form of weighted Opial inequalities for linear differential operator with related extreme cases. Later on Koliha and Pe\v cari\'c \cite{KP} gave the more generalized form  of the Opial inequalities established in \cite{ANAPE} with applications for fractional derivatives. Our purpose is to give applications of weighted Opial inequalities for Widder derivatives. Also we will give the connections of our results in this paper with \cite{ANAPE} and \cite{ANPA}(see also \cite[Chapter 12]{ANAI}) and we will show how the results in this paper generalizes the results in \cite{ANAPE} and \cite{ANPA}.

The following hypotheses are assumed throughout this section: Let $I$ be a closed interval in $\mathbb{R}$ and a be fixed point in $I$, let $\Phi$ be a continuous function nonnegative on $I \times I$, and let $y, h \in C(I)$. We assume that the following condition involving $\Phi,$ $h$ and $y$ is satisfied:
\begin{equation}\label{2.1}
|y(x)|\leq \left|\int\limits_{a}^{x}\Phi(x,t)|h(t)|dt\right|,\qquad x\in I.
\end{equation}

Some typical example of (\ref{2.1}) is given below.
\begin{example}
Let $K$ be a continuous function on $I\times I$ and let $y$ be defined by
$$y(x)=\int\limits_{a}^{s}K(s,t)h(t)dt, s\in I.$$
Then $(\ref{2.1})$ holds with $\Phi(s,t)=|K(s,t)|.$ A useful modification of this example--easier to attain in practice-- is obtained when a function $z\in C(I)$ defined by
$$z(s)= \int\limits_{a}^{s}K(s,t)h(t)dt$$
satisfies a inequality $|z(t)|\geq|y(t)|.$ Again $(\ref{2.1})$ holds with $\Phi(s,t)=|K(s,t)|.$
\end{example}

In \cite{ANAPE}, the results yields the Opial-type inequalities for linear differential operator (see \cite{AG}, \cite{AP2}, \cite{AN}).

\begin{example}\label{ex1}
 Let
$$L=\sum\limits_{j=0}^{n-1}a_{j}(t)D^{j}+D^{n},\quad t\in I,$$
be the linear differential operator with $a_{j}\in C(I),$ let $h\in C(I).$ Let $y_{1}(x),$...,$y_{n}(x)$ be the set of lineary independent solution of $L_{y}=0$ and here is the associated Green's function for $L$ is
$$
G(x,t):=\frac{\left|\begin{array}{ccccc}
               y_{1}(t) & \cdot & \cdot & \cdot & y_{n}(t) \\
               y_{1}^{\prime}(t) & \cdot & \cdot & \cdot & y_{n}^{\prime}(t) \\
               \cdot & \cdot &  &  &\cdot  \\
               \cdot &  & \cdot &  & \cdot \\
               \cdot &  &  & \cdot & \cdot \\
               y_{1}^{(n-2)}(t) &  &  & \cdot & y_{n^{(n-2)}}(t) \\
               y_{1}(x) & \cdot & \cdot & \cdot & y_{n}(x)
             \end{array}
\right|}{\left|\begin{array}{ccccc}
               y_{1}(t) & \cdot & \cdot & \cdot & y_{n}(t) \\
               y_{1}^{\prime}(t) & \cdot & \cdot & \cdot & y_{n}^{\prime}(t) \\
               \cdot & \cdot &  &  &\cdot  \\
               \cdot &  & \cdot &  & \cdot \\
               \cdot &  &  & \cdot & \cdot \\
               y_{1}^{(n-2)}(t) &  &  & \cdot & y_{n^{(n-2)}}(t) \\
               y_{1}(t) & \cdot & \cdot & \cdot & y_{n}(t)
             \end{array}
\right|},
$$
which is continuous function on $I^{2}.$
It is known that
$$y(x)= \int\limits_{a}^{x}G(x,t)h(t)dt$$
is the unique solution to the initial value problem
$$L_{y}=h,\quad y^{(j)}(a)=0,\quad j=0,1,...,n-1.$$
Then (\ref{2.1}) is satisfied for $y$ and $h$ with with $\Phi(s,t)=|G(s,t)|.$
\end{example}

The rest of the paper is organized in the following way: In Section 2, we give the preparatory weighted Opial inequalities for general kernels. In Section 3, we give the applications for Widder derivatives and we will show that the results in this paper generalizes the results in \cite{ANPA} and \cite{ANAPE}. In Section 4, we give some concluding remarks about the applications of fractional derivatives and linear differential operator and and will show how the results of \cite{KP} generalizes the results of \cite{AN} and \cite{ANAPE}.
\bigskip
\section{Preparatory inequalities}
The following result is proved in \cite{KP}.
\begin{theorem}\label{main}
Assume that $(\ref{2.1})$ holds. Let $x\in I,$ let $\alpha,\beta>0,$ $r>\max(1,\alpha),$ and let $u,v\in C(I)$ be such that $u(s)\geq0,\,\,v(s)>0$ for all $s\in I.$ Then
 $$\left|\int\limits_{a}^{x}u(s)|y(s)|^{\beta}|h(s)|^{\alpha}ds\right|\leq
 C(x)\left|\int\limits_{a}^{x}v(s) |h(s)|^{r}ds\right|^{\frac{\alpha+\beta}{r}},$$
where
\begin{equation}\label{C}
C(x):=\left(\frac{\alpha}{\alpha+\beta}\right)^{\frac{\alpha}{r}}
\left(\int\limits_{a}^{x}\left(u^{r}(s)v^{-\alpha}(s)\right)^{\frac{1}{r-\alpha}}
P(s)^{\frac{\beta(r-1)}{r-\alpha}}ds\right)^{\frac{r-\alpha}{r}},
\end{equation}
and
\begin{equation}\label{P}
P(s):=\left|\int\limits_{a}^{s}v(t)^{-\frac{1}{r-1}} \Phi(s,t)^{\frac{r}{r-1}}dt\right|.
\end{equation}
\end{theorem}

Particularly in Theorem \ref{main} for $r=2,$ the following specialization is obtained (see also \cite[Corollary 2]{ANAPE}).
\begin{corollary}\label{corav}
Assume that $(\ref{2.1})$ holds. Let $x\in I,$ let $\beta>0,$ $0<\alpha<2.$ Let $u,v\in C(I)$ be such that $u(s)\geq0,\,\,v(s)>0$ for all $s\in I.$ Then
 $$\left|\int\limits_{a}^{x}u(s)|y(s)|^{\beta}|h(s)|^{\alpha}ds\right|\leq
 \widetilde{C}(x)\left|\int\limits_{a}^{x}v(s) |h(s)|^{2}ds\right|^{\frac{\alpha+\beta}{2}},$$
where
\begin{equation}\label{C1}
\widetilde{C}(x):=\left(\frac{\alpha}{\alpha+\beta}\right)^{\frac{\alpha}{2}}
\left(\int\limits_{a}^{x}\left(u^{2}(s)v^{-\alpha}(s)\right)^{\frac{1}{2-\alpha}}
\widetilde{P}(s)^{\frac{\beta}{2-\alpha}}ds\right)^{\frac{2-\alpha}{2}}
\end{equation}
and
\begin{equation}\label{P1}
\widetilde{P}(s):=\left|\int\limits_{a}^{s}v(t)^{-1} \Phi(s,t)^{2}dt\right|.
\end{equation}
\end{corollary}

The following extreme case analogue to the \cite[Proposition 1]{ANAPE} is given here.
\begin{theorem}\label{2.3}
Assume that $(\ref{2.1})$ holds. Let $x\in I, \alpha,\beta>0$ and $r>\max(1,\alpha),$ and let $u,v\in C(I)$ be such that $u(s)\geq0, v(s)>0$ for all $s\in I.$ Then
$$\left|\int\limits_{a}^{x}u(s)|y(s)|^{\beta}|h(s)|^{\alpha}ds\right|\leq
\int\limits_{a}^{x}u(w)\left|\int\limits_{a}^{x}v(t)k(w,t)dt \right|^{\frac{r-\alpha}{r}}dw\,\| v\|_{\infty}^{\beta}\cdot\| h\|_{\infty}^{\alpha+\beta}.$$
\end{theorem}

Following \cite{ANAPE}, we consider a situation when the exponents $\alpha,$ $\beta$ and $r$ in Theorem \ref{main} are not necessarily positive. In this case the inequality (\ref{2.1}) must be strengthened to equality.
\begin{equation}\label{1aeq}
|y(s)|= \left| \int\limits_{a}^{x}\Phi(s,t)|h(t)|dt\right|,
\end{equation}
where $k$ is again a non-negative continuous function on $I\times I$ and $y,h\in C(I).$ As before, $a$ is fixed point in the interval $I.$

\begin{theorem}\label{2.4}
Assume that (\ref{1aeq}) holds. let $x\in I$, let $u,v\in C(I)$ be such that $u(s)\geq0,v(s)>0$ for all $s\in I.$ Let $C(x)$ be defined by (\ref{C}) and (\ref{P}). Consider real number $\alpha,\beta,r $ and the following relations:\\

${\rm (i)}$  $r>1,\beta>0, 0<\alpha<r;$

{\rm(ii)} $r<\alpha<0, \beta<0;$

{\rm(iii)}$-\alpha<\beta<0, 0<\alpha<r<1;$

{\rm(iv)} $\beta>0,0<r<\min(\alpha,1);$

{\rm(v)} $\alpha<0<r<1,0<\beta<-\alpha;$

{\rm(vi)} $\beta<0, \alpha<0, r>1;$

{\rm(vii)} $1<r<\alpha, -\alpha<\beta<0;$

{\rm(viii)} $\beta>0, r<0<\alpha;$

{\rm(ix)} $\alpha<r<0, 0<\beta<-\alpha.$\\

\noindent If one of the condition ${\rm(i)-(iii)}$ satisfied, then
\begin{equation}
\left|\int\limits_{a}^{x}u(s)|y(s)|^{\beta}|h(s)|^{\alpha}ds\right|\leq
C(x)\left|\int\limits_{a}^{x}V(s) |h(s)|^{r}ds\right|^{\frac{\alpha+\beta}{r}}.
\end{equation}
If one of the conditions ${\rm(iv)-(ix)}$ is satisfied, then
\begin{equation}
\left|\int\limits_{a}^{x}u(s)|y(s)|^{\beta}|h(s)|^{\alpha}ds\right|\geq
C(x)\left|\int\limits_{a}^{x}V(s) |h(s)|^{r}ds\right|^{\frac{\alpha+\beta}{r}}.
\end{equation}
\end{theorem}

\bigskip

\section{Main Results}
First it is necessary to give some important details about Widder derivatives (see \cite{WIDD}).\\
Let $f,u_{0},u_{1},...,u_{n}\in C^{n+1}([a,b]),n\geq0,$ and the Wronskians

$$
W_{i}(x):=W[u_{0}(x),u_{1}(x),...,u_{i}(x)]=\left|\begin{array}{ccccc}
               u_{0}(x) & \cdot & \cdot & \cdot & u_{i}(x) \\
               u_{o}^{\prime}(x) & \cdot & \cdot & \cdot & u_{i}^{\prime}(x) \\
               \cdot & \cdot &  &  &\cdot  \\
               \cdot &  & \cdot &  & \cdot \\
               \cdot &  &  & \cdot & \cdot \\
               u_{0}^{(i)}(x) & \cdot & \cdot & \cdot & u_{i}^{(i)}(x)
             \end{array}
\right|,$$
$i=0,1,...,n.$ Here $W_{0}(x)=u_{0}(x).$ Assume $W_{i}(x)>0$ over $[a,b], i=0,1,...,n.$ For $i\geq0,$ the differential operator of order $i$ (Widder derivative):
$$L_{i}f(x):=\frac{W[u_{0}(x),u_{1}(x),...,u_{i-1}(x),f(x)]}{W_{i-1}(x)},$$
$i=1,...,n+1;L_{0}f(x)=f(x)$ for all $x\in [a,b].$

Consider also
$$
g_{i}(x,t):=\frac{1}{W_{i}(t)}\left|\begin{array}{ccccc}
               u_{0}(t) & \cdot & \cdot & \cdot & u_{i}(t) \\
               u_{0}^{\prime}(t) & \cdot & \cdot & \cdot & u_{i}^{\prime}(t) \\
               \cdot & \cdot &  &  &\cdot  \\
               \cdot &  & \cdot &  & \cdot \\
               \cdot &  &  & \cdot & \cdot \\
               u_{0}(x) & \cdot & \cdot & \cdot & u_{i}(x)
             \end{array}
\right|,
$$
$i=1,2,...,n;$ $g_{0}(x,t):=\frac{u_{0}(x)}{u_{0}(t)}$ for all $x,t\in [a,b].$
We also mention the generalized Widder-Talylor's formula, see \cite{WIDD}(see also \cite{ANAI}).

\begin{theorem}
Let the functions $f, u_{0},u_{1},...,u_{n}\in C^{n+1}([a,b]),$ and the Wronkians $W_{0}(x),W_{1}(x),...,W_{n}(x)>0$ on $[a,b], x\in[a,b].$ Then for $t\in [a,b]$ we have
$$f(x)=f(t)\frac{u_{0}(x)}{u_{0}(t)}+L_{1}f(t)g_{1}(x,t)+...+L_{n}f(t)g_{n}(x,t)+R_{n}(x)
$$
where
$$
R_{n}(x):=\int\limits_{t}^{x}g_{n}(x,t)L_{n+1}f(s)ds.$$
\end{theorem}

For example (see \cite{WIDD}) one could take $u_{0}(x)=c>0.$ If $u_{i}(x)=x^{i},i=0,1,...,n,$ defined on [a,b], then
$$L_{i}f(t)=f^{(i)}(t)\,\,\,\,and\,\,\,\, g_{i}(x,t)=\frac{(x-t)^{i}}{i!},\quad t\in [a,b].$$

We need

\begin{corollary}  By additionally assuming for fixed $x_{0}\in[a,b]$ that $L_{i}f(x_{0})=0,i=0,1,...,n,$ we get that
\begin{equation}
f(x):=\int\limits_{x_{0}}^{x}g_{n}(x,t)L_{n+1}f(t)dt\quad for \,\,all \,\,x\in[a,b].
\end{equation}
\end{corollary}

\begin{theorem}\label{t1}
Let $s\in I,$ let $\alpha,\beta>0,$ $r>\max(1,\alpha),$ and let $u,v\in C(I)$ be such that $u(s)\geq0,\,\,v(s)>0$ for all $s\in I.$ Then
\begin{equation}\label{r1}
\left|\int\limits_{a}^{x}u(s)|f(s)|^{\beta}|L_{n+1}f(s)|^{\alpha}ds\right|
\leq C_{1}(x)\left|\int\limits_{a}^{x}v(s) |L_{n+1}f(s)|^{r}ds\right|^{\frac{\alpha+\beta}{r}},
\end{equation}
where
\begin{equation}\label{C2}
C_1(x):=\left(\frac{\alpha}{\alpha+\beta}\right)^{\frac{\alpha}{r}}
\left(\int\limits_{a}^{x}(u^{r}(s)v^{-\alpha}(s))^{\frac{1}{r-\alpha}}
P_{1}(s)^{\frac{\beta(r-1)}{r-\alpha}}ds\right)^{\frac{r-\alpha}{r}}
\end{equation}
and
\begin{equation}\label{P2}
P_1(s):=\left|\int\limits_{a}^{s}\left(v(t)^{-\frac{1}{r-1}} g_{n}(s,t)\right)^{\frac{r}{r-1}}dt\right|.
\end{equation}
\end{theorem}

\proof Applying Theorem \ref{main} with $\Phi(s,t)=g_{n}(s,t),$ $y=f$ and $h=L_{n+1}f,$ we get the inequality $(\ref{r1}).$ \qed

\begin{remark}
\rm If we take $\alpha=\beta=1$ and the weight functions $u=v=1$ and $r=q$ in Theorem \ref{t1}, we get the exactly Corollary 12.8 of \cite{ANAI}.
So the Theorem \ref{t1} is the generalized version of the Corollary 12.8 of \cite{ANAI}. Also it is observed that the Corollary 12.8 of \cite{ANAI} is the special case of Theorem 2.4 of \cite{KP}.
\end{remark}

In particular in Theorem \ref{t1} if we take $r=2$ we have the following result.
\begin{corollary}\label{cora3}
 Let $s\in I,$ let $\beta>0, 0<\alpha<2.$ Let $u,v\in C(I)$ be such that $u(s)\geq0,\,\,v(s)>0$ for all $s\in I.$ Then
 $$\left|\int\limits_{a}^{x}u(s)|f(s)|^{\beta}|L_{n+1}f(s)|^{\alpha}ds\right|\leq
 \widetilde{C_{1}}(x)\left|\int\limits_{a}^{x}v(s) |L_{n+1}f(s)|^{2}ds\right|^{\frac{\alpha+\beta}{2}},$$
where
\begin{equation}\label{C13}
\widetilde{C_{1}}(x):=\left(\frac{\alpha}{\alpha+\beta}\right)^{\frac{\alpha}{2}}
\left(\int\limits_{a}^{x}(u^{2}(s)v^{-\alpha}(s))^{\frac{1}{2-\alpha}}
\widetilde{P}_{1}(s)^{\frac{\beta}{2-\alpha}}ds\right)^{\frac{2-\alpha}{2}}
\end{equation}
and
\begin{equation}\label{P13}
\widetilde{P_{1}}(s):=\left|\int\limits_{a}^{s}v(t)^{-1} g_{n}(s,t)^{2}dt\right|.
\end{equation}
\end{corollary}

\begin{remark}
\rm If we take $\alpha=\beta=1$ and the weight functions $u=v=1$ in Corollary \ref{cora3}, we get the Corollary 12.9 of \cite{ANAI}.
\end{remark}

\begin{theorem}\label{thext}
Let $s\in I, \alpha,\beta>0$ and $r>\max(1,\alpha),$ and let $u,v\in C(I)$ be such that $u(s)\geq0, v(s)>0$ for all $s\in I.$ Then
\begin{multline}
\left|\int\limits_{a}^{x}u(s)|f(s)|^{\beta}|L_{n+1}f(s)|^{\alpha}ds\right|\\
\leq\int\limits_{a}^{x}u(w)\left|\int\limits_{a}^{x}v(t)g_{n}(w,t)dt \right|^{\frac{r-\alpha}{r}}dw\,\| v\|_{\infty}^{\beta}\cdot\| L_{n+1}f(s)\|_{\infty}^{\alpha+\beta}.
\end{multline}
\end{theorem}

\proof The proof is similar to proof of Theorem \ref{t1}.\qed

\begin{remark}
\rm If we take $\alpha=\beta=1$ and the weight functions $u=v=1$ and $r=2$ in Theorem \ref{thext}, we get the exactly Corollary 12.10 of \cite{ANAI}.
So the Theorem \ref{thext} is the generalized version of the Corollary 12.10 of \cite{ANAI}.
\end{remark}

\begin{theorem}\label{7c}
Let $s\in I$, let $u,v\in C(I)$ be such that $u(s)\geq0,v(s)>0$ for all $s\in I.$ Let $C(x)$ be defined by (\ref{C}) and (\ref{P}). Consider real number $\alpha,\beta,r $ and the following relations:\\

${\rm (i)}$  $r>1,\beta>0, 0<\alpha<r;$

{\rm(ii)} $r<\alpha<0, \beta<0;$

{\rm(iii)}$-\alpha<\beta<0, 0<\alpha<r<1;$

{\rm(iv)} $\beta>0,0<r<\min(\alpha,1);$

{\rm(v)} $\alpha<0<r<1,0<\beta<-\alpha;$

{\rm(vi)} $\beta<0, \alpha<0, r>1;$

{\rm(vii)} $1<r<\alpha, -\alpha<\beta<0;$

{\rm(viii)} $\beta>0, r<0<\alpha;$

{\rm(ix)} $\alpha<r<0, 0<\beta<-\alpha.$\\

\noindent If one of the condition ${\rm(i)-(iii)}$ satisfied, then
\begin{equation}
\left|\int\limits_{a}^{x}u(s)|f(s)|^{\beta}|L_{n+1}f(s)|^{\alpha}ds\right|\leq
C(x)\left|\int\limits_{a}^{x}v(s) |L_{n+1}f(s)|^{r}ds\right|^{\frac{\alpha+\beta}{r}}.
\end{equation}
If one of the conditions ${\rm(iv)-(ix)}$ is satisfied, then
\begin{equation}
\left|\int\limits_{a}^{x}u(s)|f(s)|^{\beta}|L_{n+1}f(s)|^{\alpha}ds\right|\geq
C(x)\left|\int\limits_{a}^{x}v(s) |L_{n+1}f(s)|^{r}ds\right|^{\frac{\alpha+\beta}{r}}.
\end{equation}
\end{theorem}

\proof The proof is similar to proof of Theorem \ref{t1}.\qed

\bigskip

\section{Concluding Remarks}
If we replace the general kernels $\Phi(s,t)$ by particular Green's function $G(s,t)$ associated to linear differential operator $L$ in Theorem \ref{main}, Corollary \ref{corav}, Theorem \ref{2.3} and Theorem \ref{2.4} we get the \cite[Corollary 1]{ANAPE}, \cite[Corollary 2]{ANAPE}, \cite[Corollary 3]{ANAPE} and \cite[Theorem 4]{ANAPE}(also \cite[Theorem 5]{ANAPE}) respectively. With the same replacement additionally if we take $\alpha=\beta=1$ and the weight functions $u=v=1$ and $r=q$ in Theorem \ref{main}, Corollary \ref{corav},  we can get the exactly  Corollary 13.5 and Corollary 13.5 of \cite{AN}. So we can say that the results of \cite{KP} cover much more general situation and generalizes the results of \cite{AN}, \cite{ANPA} and \cite{ANAPE}.

Moreover talking about applications, in \cite{KP}, the applications for Riemann-Liouville fractional derivatives are given. Such type of applications can be given for Canavati fractional derivatives and Caputo fractional derivatives but here we omit the details. As a special case for particular value of $\alpha$ and $\beta$ and weight functions, we can get the corresponding applications for \cite{ANPA}.

\end{document}